\documentclass{jcdOF}
\usepackage{amsfonts,amssymb,amsmath}
  \usepackage{paralist}
  \usepackage[misc]{ifsym}
  \usepackage{graphics} %% add this and next lines if pictures should be in esp format
  \usepackage{epsfig} %For pictures: screened artwork should be set up with an 85 or 100 line screen
\usepackage{graphicx}  \usepackage{epstopdf}%This is to transfer  figure to .pdf figure; please compile your paper using PDFLeTex or PDFTeXify.
 \usepackage[colorlinks=true]{hyperref}
\hypersetup{urlcolor=blue, citecolor=red}

\def\currentvolume{X}
 \def\currentissue{X}
  \def\currentyear{200X}
   \def\currentmonth{XX}
    \def\ppages{X--XX}
     \def\DOI{10.3934/jcd.2022017}

\newtheorem{theorem}{Theorem}[section]

\theoremstyle{definition}
\newtheorem{definition}[theorem]{Definition}

\title[Invariant sets and stability of limit cycles] %Use the shortened version of the full title
      {Finding positively invariant sets and proving exponential stability of limit cycles using Sum-of-Squares decompositions}

\author[Elias August and Mauricio Barahona]{}


\thanks{$^*$ Corresponding author: Elias August, eliasaugust@ru.is}

\begin{document}
\maketitle
\centerline{\scshape
Elias August$^{{\href{mailto:eliasaugust@ru.is}{\textrm{\Letter}}}*1}$
and Mauricio Barahona$^{{\href{mailto:m.barahona@imperial.ac.uk}{\textrm{\Letter}}}2}$}
% Enter the first author's name and address:
%\centerline{\scshape Elias August$^*$}
\medskip
{\footnotesize
% please put the address of the first author
 \centerline{$^1$Department of Engineering, Reykjavik University, Menntavegur 1, 102 Reykjavik, Iceland}
  % \centerline{}
  % \centerline{}
} % Do not forget to end the {\footnotesize by the sign }

\medskip

%\centerline{\scshape Mauricio Barahona}
\medskip
{\footnotesize
 % please put the address of the second  and third author
 \centerline{$^2$Department of Mathematics, Imperial College London,}
   \centerline{South Kensington Campus, London SW7 2AZ, United Kingdom}
  % \centerline{}
}

\bigskip

% The name of the associate editor will be entered by an editorial staff
% "Communicated by the associate editor name" is not needed for special issue.
% \centerline{(Communicated by the associate editor name)}

%The abstract of your paper
\begin{abstract}
The dynamics of many systems from physics, economics, chemistry, and biology can be modelled through polynomial functions. In this paper, we provide a computational means to find positively invariant sets of polynomial dynamical systems by using semidefinite programming to solve \emph{sum-of-squares (SOS)} programmes.  With the emergence of SOS programmes, it is possible to efficiently search for Lyapunov functions that guarantee stability of polynomial systems. Yet, SOS computations often fail to find functions, such that the conditions hold in the entire state space. We show here that restricting the SOS optimisation to specific domains enables us to obtain positively invariant sets, thus facilitating the analysis of the dynamics by considering separately each positively invariant set. In addition, we go beyond classical Lyapunov stability analysis and use SOS decompositions to computationally implement sufficient positivity conditions that guarantee existence, uniqueness, and exponential stability of a limit cycle. Importantly, this approach is applicable to systems of any dimension and, thus, goes beyond classical methods that are restricted to two-dimensional phase space. We illustrate our different results with applications to classical systems, such as the van der Pol oscillator, the Fitzhugh-Nagumo neuronal equation, and the Lorenz system.
\end{abstract}

\section{Introduction}

For the understanding of natural dynamical systems, mathematical modelling is of utmost importance. The dynamics of systems from physics (e.g., the van der Pol oscillator and Lorenz equations to cite but two), chemistry (e.g., chemical reaction networks assuming mass action kinetics~\cite{Feinberg06}), neurology (e.g., the FitzHugh-Nagumo model), epidemiology (e.g., susceptible-infected-recovered (SIR) epidemic models) or ecology (e.g., Lotka-Volterra predator-prey models) can be modelled through polynomial functions.  Such functions are also used in many other fields including economics, pharmacokinetics, transportation and communication. In this paper, we provide a computational means to find positively invariant sets of such
dynamical systems by using \emph{sum-of-squares (SOS)} programmes. Finding positively invariant sets of a dynamical system, if they
exist, is important because they decompose the state space into subsets where solution trajectories are trapped. Therefore, one can simplify the
analysis of the system by considering the dynamics in each positively invariant set. Information about the shape of the positively invariant sets is also significant, for example, to provide necessary conditions for switching behaviour of a bistable system or to guarantee converging behaviour.

An important tool for the stability analysis of dynamical systems is Lyapunov stability theory. With the emergence of SOS programmes, it is possible to efficiently search for Lyapunov functions for polynomial dynamical systems~\cite{Pablo00}. SOS programmes provide sufficient positivity conditions that guarantee stability. However, SOS computations often fail to show that the conditions hold over the entire state space, and therefore restricting the analysis to smaller domains can be of great benefit to characterise the dynamics of the system. We consider here this issue extensively through the computation of attracting sets, repelling sets, and positively invariant sets of different systems of interest, which are used for the characterisation of their dynamics.
As a second contribution, in this paper, we go beyond classical Lyapunov stability analysis and establish SOS programmes to search for a matrix function that guarantees existence, uniqueness, and exponential stability of a limit cycle. The SOS programmes presented here are based on a sufficient positivity condition established by Peter Giesl~\cite{Giesl03,Giesl204}, and were already developed and applied computationally in the PhD thesis of the first author~\cite{August07}. In~\cite{MANCHESTER201432}, essentially the same condition (termed there the transverse contracting condition) was derived, examined theoretically, and applied to an example to show computationally the existence of an exponentially stable limit cycle for a system of dimension two. Here, we present a systematic approach to use SOS programmes to establish the presence of stable limit cycles. Importantly, by finding the positively invariant set where the limit cycle lies, and by requiring conditions to hold only in this set thereby easing the computation, we are able to prove the stability of limit cycles for systems of higher dimension.  Indeed, SOS programmes have recently been used to establish region of attractions but have been demonstrated using a system of dimension two only~\cite{8959232}. On the other hand, the method presented recently in~\cite{SCHLOSSER2021109900} does offer means to obtain excellent approximation for attractors of higher dimensions using semidefinite programming; hence, we think of the approach presented in this paper as an alternative.
Thus, the main difference of the research presented here to previous work is the presentation of a concise framework for establishing bounding regions and their use to determine stability of limit cycles, as well as a `proof of concept' of its applicability to systems of dimension larger than 2.

The structure of the paper is as follows. First, we provide some brief mathematical background on the computational tools that we use (Section~\ref{semi}). Section~\ref{finding} deals with finding positively invariant sets of dynamical systems. Section~\ref{exponential} shows how to computationally establish the existence of an exponentially stable limit cycle. We illustrate our results through applications to the van der Pol oscillator, the Lorenz system, and the FitzHugh-Nagumo neocortical model~\cite{Wilson99}. Finally, Section~\ref{conc} concludes the paper.

\section{Semidefinite programming and the sum of squares decomposition}\label{semi}
The main computational tool used in this paper is optimisation through semidefinite programmes. Programmes of this type can be solved efficiently using interior-point methods. (The interested reader is referred to reference~\cite{LVSB96} and the excellent textbook by the same authors~\cite{ConOpti}.) In semidefinite programming, we replace the nonnegative
orthant constraint of linear programming by the cone of positive
semidefinite matrices ($S \geq 0$) and pose the following minimisation problem:
\begin{align}\label{may11-1}
&\text{minimise} \quad  c^T x\nonumber\\
&\text{subject to} \quad  S(x)\geq0\ \nonumber\\
&\text{where} \quad S(x)=S_0+\sum^n_{i=1}x_iS_i
\end{align}
Here, $x\in\mathbb{R}^n$ is the free variable, and the so called \emph{problem
data}, which are given, are the vector $c\in \mathbb{R}^n$ and the
matrices $S_j\in\mathbb{R}^{m\times m}$, $j=0,\ldots,n$. Note that convexity of the set of symmetric positive semidefinite matrices in (\ref{may11-1}) implies that the minimisation problem has a global minimum.

\vspace{-0.1cm}
\subsection{Sum of squares decomposition}\label{sosd}
In problems dealing with polynomials, sometimes we are required to test for positivity. It is well known that testing positivity of a polynomial is NP-hard~\cite{KGMSNK87,PAP03}, except in very particular cases. However, the requirement of positivity can be relaxed to the condition that the polynomial function is a SOS. Clearly, this is only a \emph{sufficient} condition for positivity, i.e., a function can be positive without being a SOS, so that if a SOS is not found, we cannot make any definite statement about positivity. Although SOS conditions can be, at times, conservative, they
allow the use of computationally efficient methods to search the space of polynomials
for those polynomials that fulfil the SOS conditions.

Let  $F(x)$ be a real-valued polynomial function of degree $2d$ with $x\in\mathbb{R}^n$. A sufficient condition for $F(x)$ to be
nonnegative is that it can be decomposed into a SOS~\cite{PAP03}:
\begin{equation*}\label{may10-4}
F(x)=\sum_i f_i^2(x)\geq0
\end{equation*}
where $f_i(x)$ are polynomial functions. Importantly, SOS decompositions have an algebraic characterisation~\cite{Pablo00,PAP03}: $F(x)$ is a SOS if and only if there exists a positive semidefinite matrix $R \geq 0$ such that
\begin{equation*}
F(x)={\chi}^T  R \, \chi,\ \chi=[1,\ x_1,\ x_2,\ldots,\ x_n, x_1x_2,\ldots,\ x^d_n],
\end{equation*}
where the vector of monomials $\chi$ has length
$\ell=\binom{n+d}{d}$. Note that $R$ is not necessarily unique.
Furthermore, the equality $\sum_i f_i^2(x)=\chi^\mathrm{T}R\chi$ imposes certain constraints on $R$ of the form trace$(A_jR)=c_j$, where $A_j$ and $c_j$ are appropriate matrices and constants, respectively. For an illustration of these constraints, see Example 3.5 in~\cite{PAP03}.

To find $R$, we then solve the optimisation problem associated with the following semidefinite programme,
which is the dual of the one given by (\ref{may11-1}):
\begin{eqnarray}\label{sos-dec}
\mathrm{minimise}\ &&\mathrm{trace}(A_0R) \nonumber\\
\mathrm{subject}\ \mathrm{to}\ &&\mathrm{trace}(A_jR)=c_j,\ j=1,\ldots,m \nonumber\\
&& R\geq0
\end{eqnarray}

The following additional property will be important throughout our work below, since it allows us to make statements about positivity over specific regions~\cite{henrion_sostools_2005}:
\begin{align*}
&\text{If} \quad F(x)+p(x)h(x) =\sum_i g_i^2(x)\geq0, \\
& \qquad \text{where} \quad p(x) \geq 0 \quad
\text{and} \quad h(x) =\left\{\begin{array}{l}\leq0 \, \text{ if }\ a_i\leq x_i\leq b_i\ \forall i \\
>0 \, \text{ otherwise,}
\end{array}\right. \\
&\text{then} \quad  F(x)\geq 0  \quad \text{if }   a_i\leq x_i\leq b_i, \forall i.
\end{align*}
Here $a_i,b_i$ are constants defining intervals for each coordinate $x_i$. This property will be used to show that $F(x)$ is nonnegative in a specific region of the state or/and parameter space. Throughout this paper, we solve SOS programmes using \verb"SOSTOOLS"~\cite{PPP02}, a free, third-party \verb"MATLAB" toolbox that relies on the solver \verb"SeDuMi"~\cite{Stu99}.

\section{Finding positively invariant sets}\label{finding}
When dealing with dynamical systems that have multiple equilibria, limit cycles, chaotic attractors, and combinations of these, it is
often of interest to find \emph{positively invariant sets} for solution trajectories, provided they exist. (The following definitions are
from \cite{HKK00}.)

Let $x(t)$ be a solution of the dynamical system
\begin{equation}\label{July26-1}
\dot x(t)=f(x(t)),
\end{equation}
which we assume to exist for all $t$. Then $p$ is said to be a {\it positive limit point} of $x(t)$ if there is
a sequence $\{t_n\}$, with $t_n\rightarrow\infty$ as $n\rightarrow\infty$, such that $x(t_n)\rightarrow p$ as
$n\rightarrow\infty$. Furthermore, the set of all positive limit points of $x(t)$ is called a {\it positive limit set}.

A set $M$ is {\it invariant} with respect to $\dot{x}=f(x)$ if
\begin{equation*}
x(0)\in M \Rightarrow x(t)\in M,\ \mathrm{for\ all}\ t\in\mathbb{R}
\end{equation*}
This means, that if a solution belongs to $M$ at some time instant,
then it belongs to $M$ for all time.

A set $M$ is {\it positively invariant} with respect to $\dot{x}=f(x)$  if
\begin{equation*}
x(0)\in M \Rightarrow x(t)\in M,\ \mathrm{for\ all}\ t\geq0
\end{equation*}
Note that, in this paper, we are only interested in positively invariant sets.

In the following, we present sufficient conditions for the existence of positively invariant sets of a dynamical system~(\ref{July26-1}) where the vector field $f(\cdot)$ is polynomial, and show how the positively invariant set can be found through the solution of SOS semidefinite programmes using \verb"SOSTOOLS". (A related procedure was presented in~\cite{Antonis05}.)

\subsection{Finding positively invariant sets using SOS decompositions}
Let us consider a dynamical system given by
\begin{equation}\label{Dec18Schaffhausen}
\dot x=f(x),\ f(0)=0,\ x\in \mathcal{D},\ \mathcal{D}\subset\mathbb{R}^n,
\end{equation}
where the function $f:\mathcal{D}\rightarrow\mathbb{R}^n$ is continuously differentiable and the domain $\mathcal{D}$ contains the origin $x=0$.
In order to prove local asymptotic stability of the origin (needless to say that the result holds for any equilibrium point through a transformation of coordinates),
one is almost always required to show that there exists a continuously differentiable function $V(x)$ that is radially
unbounded, a so called \emph{Lyapunov function}, such that~\cite{HKK00}:
\begin{equation}\label{Dec18Singen}
V(0)=0,\ \dot V(0)=0,\quad  V(x)>0,\ \dot V(x)<0, \, x\in\mathcal{D}\setminus\{0\}
\end{equation}
However, a general methodology to construct such a function does not exist. Typically, if the dynamical system is a model of a real system, one uses some physical, chemical or biological insight  to propose a {Lyapunov function candidate}, and then checks whether it fulfils the conditions~\eqref{Dec18Singen}. However, when the vector field $f$ is polynomial, we can use semidefinite programming techniques to efficiently search for Lyapunov functions that are SOS polynomials (thus nonnegative) while satisfying certain constraints, e.g., the conditions in~\eqref{Dec18Singen}~\cite{Antonis05,PAP03,Pablo00}. Note that, in the following, whenever we require that function $V(x)$ is radially unbounded then we enforce the condition by requiring that $V(x)\geq\delta\|x\|_2^2$, where $\delta$ is a small positive real constant.

In the following, we use these ideas to search for a positively invariant set of \eqref{Dec18Schaffhausen} of possibly large size (with respect to some measure) around the origin
such that all solutions whose trajectories start in the interior of
the set and not on its boundary, remain in the set.
For polynomial systems, these conditions can be relaxed using SOS
and implemented computationally, as shown in the following result.
\begin{theorem}\label{thmApr2}
Let us consider a dynamical system~\eqref{Dec18Schaffhausen} with polynomial vector field $f(x)$.
If we find a radially unbounded SOS polynomial $V_a(x)$ such that $V_a(x)>0$ if $x\neq0$ and $V_a(0)=0$, a (desirably large) constant $\beta_a > 0$, and a SOS polynomial $r_1(x)$ such that
\begin{equation}\label{Sep2-eq1}
-f(x)\cdot\nabla V_a(x)+ r_1(x) \left(\|x\|_2^2-\beta_a \right )\ \mathrm{is\ SOS,}
\end{equation}
as well as a constant $\gamma_a > 0$ and a (not necessarily SOS) polynomial $r_2(x)$ such that
\begin{equation}\label{Sep2-eq2}
(V_a(x)-\gamma_a)+r_2(x) \left(\|x\|_2^2-\beta_a \right)\ \mathrm{is\ SOS},
\end{equation}
then $\mathcal{D}=\{x\in \mathbb{R}^n\ |\ V_a(x)<\gamma_a \cap \|x\|_2^2\leq\beta_a\}$ is a positively invariant set of~\eqref{Dec18Schaffhausen}.
\end{theorem}
\begin{proof}
First, note that if (\ref{Sep2-eq1}) holds and $\|x\|^2_2\leq\beta_a$
then  $\dot V_a(x)=f(x)\cdot\nabla V_a(x)\leq0$.
Next, if also (\ref{Sep2-eq2}) holds and $V_a(x)<\gamma_a$
then $r_2(x)(\|x\|_2^2-\beta_a)>0$, which implies that $\|x\|_2^2\neq\beta_a$.
It follows that $\|x\|_2^2<\beta_a$ and, thus, that $x$ is in the interior of set $\mathcal{D}$, where $V_a(x)$
is non-increasing and, since then $V_a(x)-\gamma_a<0$, in order for (\ref{Sep2-eq2}) to hold, $\|x\|_2^2$ cannot converge to $\beta_a$. This implies that if $x$ is initially in the interior of set $\mathcal{D}$ then trajectories will be such that $\|x\|_2^2<\beta_a$. Therefore, the set $\mathcal{D}$ is positively invariant.
\end{proof}

Now we consider absorbing sets, which are a special subclass of positively invariant sets, as given in the following definition.
\begin{definition} \label{dissipativity}
A dynamical system $\dot x = f(x)$ is said to be ultimately bounded if it
possesses a bounded absorbing set.
For real and positive constants $\alpha$ and $\beta$, the bounded set $\mathcal{B}:=\{x\in\mathbb{R}^n:\|x\|\leq \beta\}$ is absorbing if for all $x(0)$ for which $\|x(0)\|\leq\alpha$, there exists $t^\star=t^\star(\alpha,\beta)\geq0$ such that $x(t)\in \mathcal{B}$ for all $t\geq t^\star$.
Equivalently, the dynamical system $\dot x = f(x)$ is said to be ultimately bounded if there exists a constant $\beta_u\geq0$ and a radially unbounded Lyapunov function $V_u(x)$ such that $V_u(x)>0$ if $x\neq0$, $V_u(0)=0$, and $\dot V_u(x)= f(x)\cdot\nabla V_u(x) <0$, if  $\|x\|_2^2\geq\beta_u$. (See, for example, \cite{IDC02,HKK00}).
\end{definition}

The next theorems provide sufficient conditions to find absorbing sets in polynomial dynamical systems based on SOS decompositions.
\begin{theorem}\label{cor1}
Let $f(x)$ be the polynomial vector field of a dynamical system~\eqref{Dec18Schaffhausen}.
If there exist constants $\beta_u \geq 0$ and $\sigma_u >0$, and SOS polynomials $p_1(x)$ and $V_u(x)$ such that
$V_u(x)$ is radially unbounded, $V_u(x)>0$ if $x\neq0$, $V_u(0)=0$, and
\begin{equation}\label{beta1}
-f(x)\cdot\nabla V_u(x) - p_1(x) \left(\|x\|_2^2-\beta_u\right)-\sigma_u\ \mathrm{is\ SOS}
\end{equation}
then the dynamical system $\dot x =f(x)$ is ultimately bounded.
\end{theorem}
\begin{proof}
If \eqref{beta1} holds, then $f(x)\cdot\nabla V_u(x)=\dot V_u(x) < 0$ if $\|x\|^2_2\geq\beta_u$, and it follows from Definition~\ref{dissipativity} that the dynamical system is ultimately bounded and it has a bounded absorbing set.
\end{proof}

Now, it follows that if \eqref{beta1} holds and there exits a
positive constant $\gamma_u$ such that $\|x\|^2_2\geq\beta_u$ if
$V_u(x)>\gamma_u$ then the set given by
$\mathcal{B}=\{x\in \mathbb{R}^n\ |\ V_u(x)\leq\gamma_u\}$ is an absorbing  set of~\eqref{Dec18Schaffhausen}. The following theorem provides sufficient conditions for the existence of such a $\gamma_u$.

\begin{theorem}\label{cor2}
Let the radially unbounded SOS polynomial $V_u(x)$ and the constant $\beta_u \geq 0$ be given such that (\ref{beta1}) holds. If there exist a constant $\gamma_u > 0 $ and a SOS polynomial $p_2(x)$ such that
\begin{equation}\label{gamma1}
-\left(V_u(x)-\gamma_u \right) + p_2(x) \left(\|x\|_2^2-\beta_u \right)\ \mathrm{is\ SOS}
\end{equation}
then $\|x\|^2_2>\beta_u$ if $V_u(x)>\gamma_u$ such that
$\mathcal{B}=\{x\in \mathbb{R}^n\ |\ V_u(x)\leq\gamma_u\}$ is an absorbing  set of~\eqref{Dec18Schaffhausen}.
\end{theorem}
\begin{proof}
If $V_u(x)>\gamma_u$ then (\ref{gamma1}) holds only if $\|x\|_2^2>\beta_u$ and $p_2(x)>0$.
\end{proof}

\begin{figure}[h]
\begin{picture}(0,680)\hspace{-9.25cm}
\includegraphics[width=1.45\textwidth]{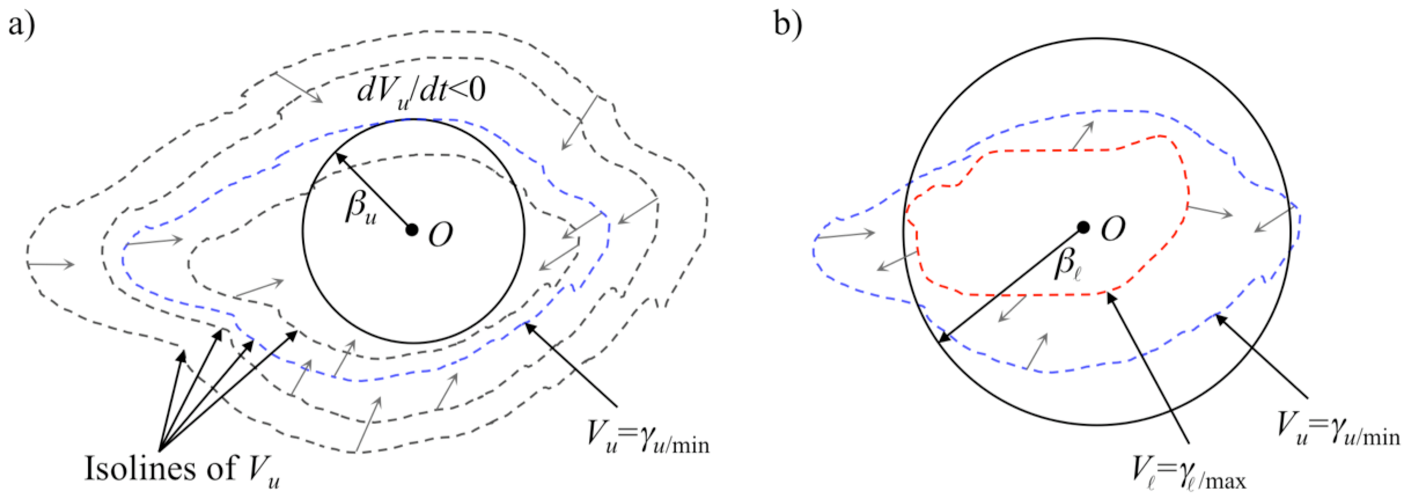}
\end{picture}
\vspace{-19cm}
\caption{\textbf{Finding positively invariant
sets.} a) The function $V_u(x)$ is constant along each of the dashed lines (\emph{isolines}). Let $V_u(x)$ and $\beta_u$ be such that (\ref{beta1}) holds. The isoline $V_u=\gamma_{u/\mathrm{min}}$ bounds the set with the smallest Euclidean distance to the origin, since $\gamma_{u/\min}$ is the minimal value of $\gamma_u$ for which (\ref{gamma1}) holds. The solution trajectories point to the interior of this set. b) Similarly, let  $V_\ell(x)$ and $\beta_\ell$ be such that (\ref{Apr2-eq1}) holds. The isoline $V_\ell=\gamma_{\ell/\mathrm{min}}$ bounds the set with the largest Euclidean distance to the origin, since $\gamma_{\ell/\max}$ is the maximal value of $\gamma_\ell$ for which (\ref{Apr2-eq2}) holds. The solution trajectories point to the exterior of this set. The region between $V_u=\gamma_{u/\mathrm{min}}$ and $V_\ell=\gamma_{\ell/\max}$ is positively invariant.}\label{isolines}
\end{figure}

Note that the constant $\beta_u$ determines the size of the absorbing set found computationally, which we aim to make as small as possible so as to provide tight bounds on the asymptotic behaviour of the system.
Therefore, computationally we solve (\ref{beta1}) repeatedly for increasingly smaller values of $\beta_u$.
Note, also, that by minimising $\gamma_u$ we decrease the Euclidean distance to the origin of the set defined by $V_u(x)-\gamma_u=0$.
Figure~\ref{isolines}a provides an illustration of the different functions in Theorem~\ref{cor1} and Theorem~\ref{cor2}.

Similarly, we can search for repelling regions, i.e., sets that will not be entered by solution trajectories.
This can be achieved by searching for a radially unbounded function $V(x)$ and a set $\mathcal{E}$ such that
\begin{equation}
V(0)=0,\ \dot V(0)=0,\quad V(x)\geq0,\ \dot V(x) \geq 0,\ x\in\mathcal{E},
\end{equation}
since this implies that, on the boundary of the set $\mathcal{E}$, the gradients of the flow $\dot V(x)=f(x)\cdot\nabla V(x)$, point towards the outside of the set or are zero. The following theorem is analogous to Theorem~\ref{thmApr2}.
\begin{theorem}\label{thmMay9}
Let us consider a dynamical system~\eqref{Dec18Schaffhausen} with polynomial vector field $f(x)$.
If we find a radially unbounded SOS polynomial $V_\ell(x)$, where $V_\ell(0)=0$, a (desirably large) constant $\beta_\ell > 0$ and a SOS polynomial $q_1(x)$ such that
\begin{equation}\label{Apr2-eq1}
f(x)\cdot\nabla V_\ell(x)+ q_1(x) \left(\|x\|_2^2-\beta_\ell \right )\ \mathrm{is\ SOS,}
\end{equation}
as well as a constant $\gamma_\ell > 0$ and a (not necessarily SOS) polynomial $q_2(x)$ such that
\begin{equation}\label{Apr2-eq2}
(V_\ell(x)-\gamma_\ell)+ q_2(x) \left(\|x\|_2^2-\beta_\ell \right)\ \mathrm{is\ SOS},
\end{equation}
then $\mathcal{E}=\{x\in \mathbb{R}^n\ |\ V_\ell(x) < \gamma_\ell \cap \|x\|_2^2\leq\beta_\ell\}$ is a repelling set of~\eqref{Dec18Schaffhausen}.
\end{theorem}
\begin{proof}
To see that this is true, first, note that if (\ref{Apr2-eq1}) holds then $\dot V_\ell(x)=f(x)\cdot\nabla V_\ell(x)\geq0$ if
$\|x\|^2_2\leq\beta_\ell$.
Next, if (\ref{Apr2-eq2}) holds and $V_\ell(x)<\gamma_\ell$ then
$q_2(x)(\|x\|_2^2-\beta_\ell)>0$, which implies that $\|x\|_2^2\neq\beta_\ell$.
It follows that $\|x\|_2^2<\beta_\ell$.
Therefore, by continuity of $V_\ell(x)$ and of the function given by $\|x\|_2^2-\beta_\ell$, if
$V_\ell(x)=\gamma_\ell$ then $\|x\|_2^2\leq\beta_\ell$, which means that
$\dot V_\ell(x)=f(x)\cdot\nabla V_\ell(x)\geq0$, as shown previously, which completes the proof.
Alternatively, Theorem~\ref{thmMay9} follows from Theorem~\ref{thmApr2} by reversing time.
\end{proof}

Similarly to before, by maximising $\gamma_\ell$ we increase the Euclidean distance to the origin of the set defined by $V_\ell(x)-\gamma_\ell=0$.
Note that the region between the isolines $V_\ell(x)=\gamma_\ell$
and $V_u(x)=\gamma_u$ is positively invariant (Figure~\ref{isolines}b). In the following, we
apply this procedure to find positively invariant sets for the Lorenz system, the FitzHugh-Nagumo model, and the van der Pol oscillator.

\subsubsection{The Lorenz system}

\begin{figure}[h]
\begin{picture}(0,220)\hspace{-5cm}
\includegraphics[width=0.75\textwidth]{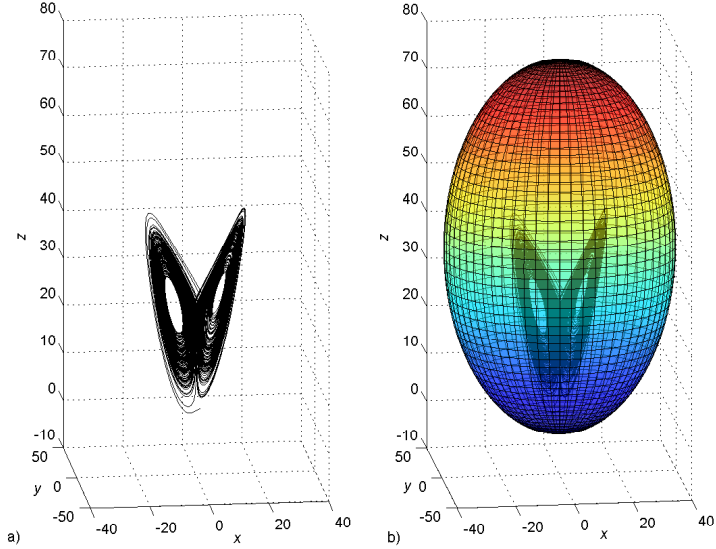}
\end{picture}%\vspace{0cm}
\caption{\textbf{The Lorenz
attractor of the system~\eqref{lor} and the classic absorbing set given by~\eqref{boundDec20}.}}
\label{LorExample}
\end{figure}

Consider the famous Lorenz system~\cite{Lorenz}, which was proposed as a simple model for circulation in the atmosphere, and has become one of the classic examples in chaos theory:
\begin{eqnarray}\label{lor}
\dot x&=& \sigma(y-x) \nonumber\\
\dot y&=& rx-y-xz\nonumber\\
\dot z&=& xy-bz
\end{eqnarray}
To study its chaotic attractor, the following system parameters are typically chosen:
$\sigma=10$, $b=8/3$, and $r=28$~\cite{Lorenz}.
A well known absorbing set for this system~\cite{Belykh04} is bounded by the ellipsoid
\begin{equation}\label{bound}
x^2+y^2+(z-\sigma-r)^2-\frac{b^2(r+\sigma)^2}{4(b-1)}=0
\end{equation}
which, for our parameters, becomes
\begin{equation}\label{boundDec20}
x^2+y^2+(z-38)^2-1540.3=0.
\end{equation}
Figure~\ref{LorExample} shows the Lorenz attractor of~\eqref{lor} and the absorbing set~\eqref{boundDec20}.

\newpage
\vspace{-3.5ex}
Let $x_1=x$, $x_2=y$, and $x_3=z-38$. We computationally implement the conditions in Theorems~\ref{cor1}~and~\ref{cor2} to find an absorbing set. First, for a given positive constant $\beta_u$, we search for functions $V_u(x)$ and $p_1(x)$ such that \eqref{beta1} holds. This feasibility problem can be implemented as follows:
\begin{eqnarray}\label{Dec14}
\mathrm{Given}\ && f(x),\ \beta_u,\ \delta>0\nonumber\\
\mathrm{search\ for}\ && V_u(x),\ p_1(x)\nonumber\\
\mathrm{subject\ to}\ && V_u(x)-\delta\|x\|^2_2\ \mathrm{is\ SOS} \nonumber\\
&&  p_1(x)\ \mathrm{is\ SOS} \nonumber\\
&& -f(x)\cdot\nabla V_u(x)-p_1(x)(\|x\|_2^2-\beta_u)\ \mathrm{is\
SOS}
\end{eqnarray}
Note that we use a (small) positive constant to guarantee that $V_u(x)>0$ if $x\neq0$. Here, we use $\delta=0.0001$.

We computationally implement programme \eqref{Dec14} in \verb"MATLAB" using the \verb"SOSTOOLS" toolbox. (We refer the reader to~\cite{PPP02} for a detailed explanation of the toolbox.)
We repeatedly solve \eqref{Dec14} for decreasing values of $\beta_u$ to obtain the smallest value of $\beta_u$ for which \eqref{beta1} holds. We find that this minimised value is $\beta_u=700$  and the corresponding function $V_u(x)$ is:
\begin{equation}\label{newLorB}
 \begin{split}
V_u(x)=&0.1324\cdot10^{-5}x_3^2x_1^4+0.5062\cdot10^{-6}x_3^4x_1^2-0.5331\cdot10^{-1}x_1x_2\\
 &+0.1282\cdot10^{-5}x_3^2x_2^2x_1^2-0.1174\cdot10^{-5}x_3^2x_2^3x_1-0.15521\cdot10^{-5}x_3^2x_2x_1^3x_1x_2\\
 &-0.5228\cdot10^{-6}x_3^4x_2x_1+0.2578\cdot10^{-5}x_1^6+0.4179\cdot10^{-6}x_2^6\\
 &+0.2961\cdot10^{-1}x_1^2x_1x_2+0.3436\cdot10^{-1}x_2^2+0.9862\cdot10^{-3}x_3^2\\
 &-0.1477\cdot10^{-6}x_1x_3+0.7294\cdot10^{-6}x_3^2x_2^4x_1x_2+0.245\cdot10^{-6}x_3^4x_2^2\\
 &+0.1964\cdot10^{-6}x_2x_3+0.152\cdot10^{-5}x_2^2x_1^4-0.5314\cdot10^{-6}x_2^5x_1x_1x_2\\
 &+0.9219\cdot10^{-6}x_2^4x_1^2-0.1534\cdot10^{-5}x_2^3x_1^3-0.8322\cdot10^{-6}x_2x_1^5
 \end{split}
\end{equation}
We also find that $p_1(x)=0.0078$; however, this function only aids the computation and we are not interested in its form.

The function $V_u(x)$ defines the shape of the absorbing set of the Lorenz system. To find this set we implement \eqref{gamma1} the following minimisation problem, where $\beta_u$ and $V_u(x)$ are the ones we determined previously:
\begin{eqnarray}\label{Dec20eq1}
\mathrm{given}\ && \beta_u,\ V_u(x)\nonumber \\
\mathrm{minimise}\ && \gamma_u\nonumber\\
\mathrm{subject\ to}\ && \gamma_u>0\nonumber \\
&& p_2(x)\ \mathrm{is\ SOS} \nonumber\\\
&&-(V_u(x)-\gamma_u)+p_2(x)(\|x\|_2^2-\beta_u)\ \mathrm{is\ SOS}
\end{eqnarray}
This programme minimises the positive constant $\gamma_u$ to `tighten' the set given by $V_u(x)=\gamma_u$. (Again, the function $p_2(x)$ only aids the computation and we are not interested in its form.) Through this procedure we finally obtain  \[\gamma_u=917.16. \]

\begin{figure}[h]
\hspace{-12cm}
\begin{picture}(0,200)
\includegraphics[width=.9\textwidth]{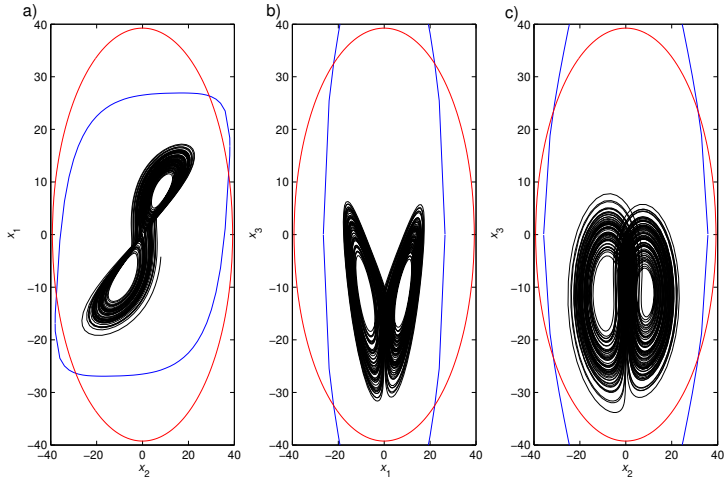}
\end{picture}\caption{\textbf{The two positively invariant sets of the Lorenz system.} The figure shows a projection of the trajectory of the Lorenz system (in black), cuts through the sphere defined by (\ref{boundDec20}) (in red) and the set given by the interior of $V_u(x)=\gamma_u$ (in blue). Note that the intersection of the sets bounded by the red and blue lines bounds the trajectories of the Lorenz system.}\label{LorBound}
\end{figure}

Figure~\ref{LorBound} shows how the set $V_u(x)\leq 917.16$
compares with \eqref{boundDec20}, by showing cuts through both absorbing sets. Note that the intersection of the two sets is a positively invariant set and, thus, the combination with our set results in a tighter bound for the trajectories of the Lorenz system. We note that (\ref{Dec14}) and (\ref{Dec20eq1}) can be combined into one optimisation problem. However, for readability, we have not done so.

\subsubsection{FitzHugh-Nagumo model}
The dynamics of neocortical neurons in humans and mammals are governed by roughly a dozen ion currents and their interplay. This leads to complex behaviour and thus, low-dimensional models have been developed in order to provide insight into the relationship between the core dynamical principles and the biophysics of the neuron~\cite{Koch99}. A popular model was proposed by the Nobel laureates Alan Lloyd Hodgkin and Andrew Huxley, and further simplified by Richard FitzHugh and Jin-Ichi Nagumo (and many others) to the following expression~\cite{MU90,Wilson99}:
\begin{eqnarray}\label{fitz}
C\dot v&=&m_\infty[v](v - E) - w + I\nonumber\\
\dot w&=&bv-\gamma w
\end{eqnarray}
where $m_\infty[v] = v(a-v)$. The first equation in \eqref{fitz} describes the changes in membrane potential $v$ which
depend on changes in membrane capacitance, ionic currents, and the stimulating current $I$ (in nA). Here the term $m_\infty[v](v - E)$ is the current contribution from ionic movement, where the constant $E$ is the steady state ion potential, and $m_\infty[v]$ is the ion activation function.
%Note that the time $t$ is measured in ms, and the voltage is measured in $10^{-2}$mV.
The parameter $a$ is neuron-dependent and leads to
different dynamical responses. The second equation in \eqref{fitz} describes the excitation and relaxation dynamics of an intrinsic variable introduced to fit experimental data from neuronal dynamics.
Below we will set $x_1=v$ and $x_2=w$ and $C=E=b=1$ and $a=4$ throughout. We will explore changes in the dynamics of the system under variation of the input current $I$ and the parameter $\gamma$.

\begin{figure}[h]\hspace{-10cm}
\begin{picture}(0,225)
\includegraphics[width=.75\textwidth]{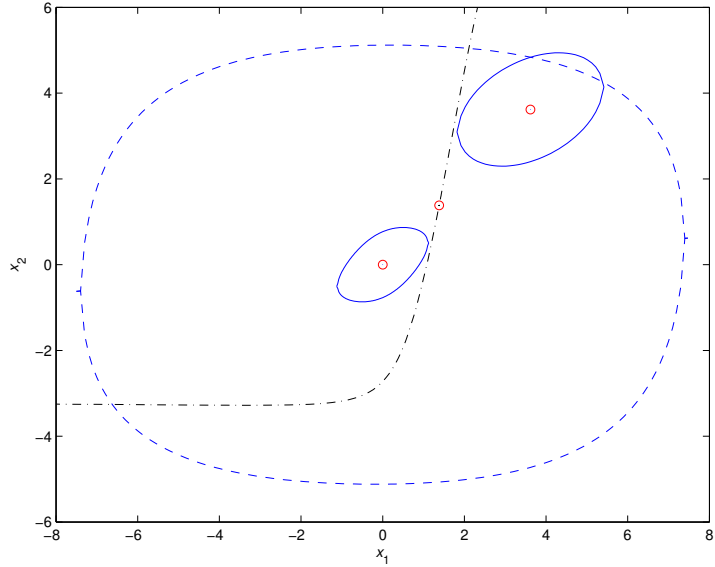}
\end{picture}\caption{\textbf{Positively invariant sets of the FitzHugh-Nagumo model~\eqref{Aug28-eq1} with $I=0$.} The three equilibrium points are indicated by red open circles. The overall attracting set~\eqref{eq:overall_attracting_setFH} is bounded by the blue dashed line. Solution trajectories around the individual equilibrium
points are bounded by solid blue lines given by~\eqref{Sep2-eq3}~and~\eqref{Sep2-eq4}.
The black dash-dotted line was calculated numerically and shows the true border between the positively invariant sets around the equilibria.
Switching between the two stable steady states
will not be observed as long as perturbations remain within the solid blue lines.
}\label{lhf}
\end{figure}

\emph{Excitable regime.}
We start by examining the FitzHugh-Nagumo system~\eqref{fitz} with the parameters above and $\gamma=1$.
Then the system~\eqref{fitz} becomes
\begin{eqnarray}\label{Aug28-eq1}
\dot x_1&=&-x_1^3+5x_1^2-4x_1- x_2+I\nonumber\\
\dot x_2&=&x_1-x_2
\end{eqnarray}
We first consider the case with no input current, $I=0$.
With these parameters, \eqref{Aug28-eq1} exhibits switching behaviour, that is, it has two stable equilibria, $(0,0)$ and $(3.618,3.618)$ and one unstable, $(1.382,1.382)$.
As above, we use \verb"SOSTOOLS" to search for an attracting set around the three equilibria that traps solution trajectories. We implement  \eqref{beta1} and \eqref{gamma1}
to find
\begin{align}
\label{eq:overall_attracting_setFH}
\beta_u&=26.2,\ \gamma_u=88.183\\
V_u(x)&=0.109x_1^2-0.089x_1x_2+0.484x_2^2+0.028x_1^4-0.007x_1^3x_2+0.045x_1^2x_2^2\nonumber\\&+0.11x_2^4
\end{align}
The boundary of the absorbing set is given by $V_u(x)=\gamma_u$ and is shown in Figure~\ref{lhf} (blue dashed line).

We then find two positively invariant sets around each of the two stable equilibria by implementing (\ref{Sep2-eq1}) and (\ref{Sep2-eq2}) in \verb"SOSTOOLS".
For the equilibrium at $(0,0)$, we get:
\begin{align}\label{Sep2-eq3}
\beta_a&=1.56,\ \gamma_a=0.381 \nonumber\\
V_a(x)&=0.422x_1^2-0.552x_1x_2+0.609x_2^2+0.066x_1^4-0.023x_1^3x_2+0.224x_1^2x_2^2\nonumber\\&+0.063x_2^4
\end{align}
To obtain the positively invariant around the
equilibrium point $(3.618,3.618)$, we shift it to the
origin through a change of variables: $y_1=x_1-3.618$ and
$y_2=x_2-3.618$, and implement (\ref{Sep2-eq1})--(\ref{Sep2-eq2}) using \verb"SOSTOOLS". This leads to the following result:
\begin{align}\label{Sep2-eq4}
\beta_a&=3.6,\ \gamma_a=0.235 \nonumber\\
V_a(y)&=0.088y_1^2-0.093y_1y_2+0.159y_2^2
\end{align}
Both bounding regions defined by $V_a = \gamma_a$ as given by (\ref{Sep2-eq3}) and
(\ref{Sep2-eq4}) are shown in Figure~\ref{lhf} (blue
solid lines).

Neural systems are excitable, that is,
if the system at equilibrium is sufficiently perturbed, the variables will display a considerable large excursion in phase space. This phenomenon can be also understood as switching between two
stable steady states: an active one (the equilibrium point corresponding to nonzero membrane potential) and and inactive one (the origin)~\cite{Koch99}.
Figure~\ref{lhf} shows that such behaviour will not be observed as long as the instantaneous perturbations remain within the solid blue circles. In other words, activation (i.e., neural ``firing'') as well as inactivation of the system are safeguarded against noise (perturbations) and require a large enough trigger (voltage build-up or drop-off) to induce the switching.

We then examine the effect of increasing the driving current $I$ in this bistable regime for the same parameters. To do this, we use Dulac's criterion:
\begin{theorem}[Dulac's criterion~\cite{JGPH83}] Let $\dot x=f(x)$ with $x\in\mathbb{R}^2$.  If there is a function $B(x)$ such that, for all $x$,
$\nabla B \cdot f$
is either positive or negative (but not both), then a periodic solution does not exist.
\end{theorem}

Using \verb"SOSTOOLS", we search for such a SOS function $B(x)$ for our choice of parameters. We find that periodic behaviour will not be observed for $1.49\leq I$. We confirmed this result numerically; also, we do not observe periodic behaviour in numerical simulations of (\ref{Aug28-eq1}) for $0\leq I<1.49$. 

\begin{figure}[h]\hspace{-11cm}
\begin{picture}(0,240)
\includegraphics[width=0.8\textwidth]{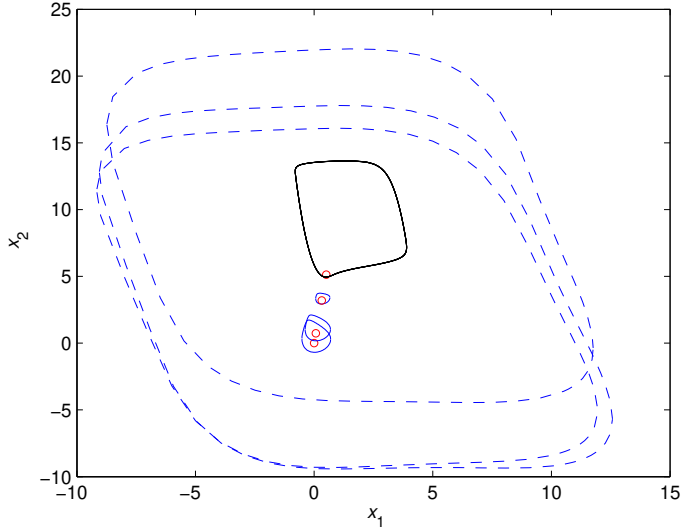}
\end{picture}\caption{\textbf{Change of the positively invariant sets of the FitzHugh-Nagumo model~\eqref{Aug28-eq1} with $I$.} The equilibrium points are presented by red circles for $I=0,1,4,6$ and the corresponding positively invariant sets around them (if they exist) are shown in blue (solid lines). The absorbing sets for $I=0,1,6$ are represented by dashed blue lines (which move up together with the corresponding equilibrium points) and the limit cycle for $I=6$ is in black.}\label{NeoCor}
\end{figure}

\emph{Periodic behaviour.}
Although periodic behaviour was not observed for the above parameters with $\gamma=1$,
periodic firing is typical of neural systems and often observed when the stimulating current $I$ is increased. Indeed, for many different sets of model parameter, the Fitzhugh-Nagumo model exhibits such behaviour. For instance, we consider the model with the same parameters above but with the parameter $\gamma=0.1$:
\begin{eqnarray}\label{Jan9-eq1}
\dot x_1&=&-x_1^3+5x_1^2-4x_1- x_2+I\nonumber\\
\dot x_2&=&x_1-0.1x_2
\end{eqnarray}
We can then use our SOS approach to provide a picture of the dynamical behaviour of (\ref{Jan9-eq1}) as we increase $I$.
In this case, for $I=0$ there exists a unique stable equilibrium point at $(0,0)$.
The fixed point remains stable and meanders for $0<I$; becomes unstable for $6\lesssim I$; and becomes stable again for $23\lesssim I$.
Using \verb"SOSTOOLS", we find positively invariant sets around the equilibrium point, given by a polynomial of degree 4, for different values of $I$. Figure~\ref{NeoCor} was created using the analysis framework presented so far in this paper.
We find that the positively invariant set shrinks as $I$ increases until it disappears completely when the equilibrium point becomes unstable (the sets for $I=0,1,4$ are shown in solid blue lines in Figure~\ref{NeoCor}). On the other hand, the absorbing set remains virtually unaffected by changes in $I$ (the sets for $I=0,1,6$ are shown in blue dashed lines).
In summary, the positively invariant regions around the fixed point in Figure~\ref{NeoCor} show that the size of the perturbation necessary to excite the system decreases as $I$ increases from $0$ to $6$. (Note that the perturbation might still not be sufficient.) The figure also shows the limits of the excursion, given by the dashed blue lines that show the bounds on the trajectories. Finally, we show that the disappearance of the positively invariant set surrounding the fixed point corresponds to the stable constant steady state ceasing to exist and becoming a stable limit cycle (stability follows since solution trajectories are bounded).
%\vspace{-1.5cm}

\subsubsection{The van der Pol
oscillator}
\label{sec:vanderPol_LC}

We have also studied the classic van der Pol oscillator, originally introduced in electronics and which is widely used to model oscillatory behaviour, e.g., of the heart~\cite{SLV04} and of neuron activity~\cite{FitzHugh}:
\begin{eqnarray}\label{van der Pol}
\dot x_1&=&x_2\nonumber\\
\dot x_2&=&-(k+x_1^2)\, x_2-x_1
\end{eqnarray}
Here, $k \in \mathbb{R}$ is a parameter. The origin $(0,0)$ is the
only fixed point of (\ref{van der Pol}): for $k>0$, the origin is
stable; for $k<0$, the origin becomes unstable and there exists a
unique stable limit cycle. Its existence and stability can be proved using Li\'enard's theorem.

We now apply our results to find a positively invariant set for (\ref{van der Pol}) for $k=-1$.
As above, we search for functions $V_u(x)$ and $p_1(x)$ such that (\ref{beta1}) holds, for increasingly smaller positive constants $\beta_u$. This feasibility problem is implemented as in (\ref{Dec14}) with $\delta=0.0001$. We find a minimised  $\beta_u=3.7$ and the corresponding function
\begin{eqnarray}\label{Aug29-eq2}
V_u(x)&=&8.361x_1^2-11.679x_1x_2+4.925x_2^2-4.098x_1^4+3.240x_1^3x_2+0.710x_1^2x_2^2\nonumber\\
&&-0.063x_1x_2^3+0.050x_2^4+0.781x_1^6+0.298x_1^4x_2^2-0.434x_1^3x_2^3+0.283x_1^2x_2^4\nonumber\\
&&-0.089x_1x_2^5+0.012x_2^6.
\end{eqnarray}

We then implement (\ref{Dec20eq1}) to tighten the set $V_u(x)=\gamma_u$ and obtain that the minimal $\gamma_u$ is given by:
\begin{equation}\label{Dec20-eq2}
\gamma_u=41.
\end{equation}
We then find a repelling set that is as large as possible by implementing~\eqref{Apr2-eq1} and~\eqref{Apr2-eq2}. We find, first, $\beta_\ell=5.8$ and
\newpage
\begin{eqnarray}\label{Aug29-eq3}
&&V_\ell(x)\nonumber\\
&=&(10.854x_1-1.942x_2-1.258x_1^3+0.517x_1^2x_2-0.244x_1x_2^2-0.105x_2^3)x_1\nonumber\\
&&+(-1.942x_1+6.094x_2+0.315x_1^3-0.571x_1^2x_2-0.159x_1x_2^2-0.447x_2^3)x_2\nonumber\\
&&+(-1.258x_1+0.315x_2+0.180x_1^3-0.093x_1^2x_2-0.030x_1x_2^2+0.017x_2^3)x_1^3\nonumber\\
&&+(0.517x_1-0.570x_2-0.092x_1^3+0.107x_1^2x_2+0.021x_1x_2^2+0.007x_2^3)x_1^2x_2\nonumber\\
&&+(-0.244x_1-0.159x_2-0.030x_1^3+0.021x_1^2x_2+0.149x_1x_2^2+0.003x_2^3)x_1x_2^2\nonumber\\
&&+(-0.105x_1-0.447x_2+0.017x_1^3+0.007x_1^2x_2+0.003x_1x_2^2+0.085x_2^3)x_2^3\nonumber\\
\end{eqnarray}
and, second, the maximal value
\begin{equation}\label{Dec20-eq3}
\gamma_\ell=12.4.
\end{equation}

\begin{figure}[h]
%\begin{picture}(0,270)\hspace{0cm}
\begin{center}
\includegraphics[width=.8\textwidth]{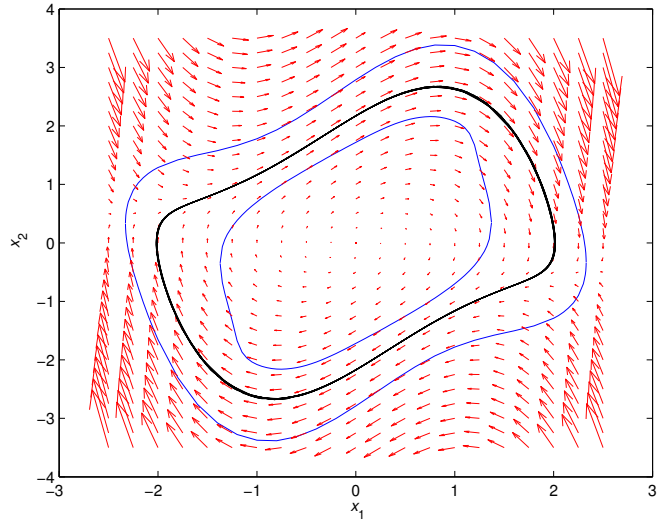}
\end{center}
%\end{picture}\vspace{-0.5cm}
\caption{\textbf{Positively invariant set and trajectories of the van der Pol oscillator.} The boundaries of the positively invariant set $\mathcal{B}$ (in blue) correspond to~\eqref{Aug29-eq2}-\eqref{Dec20-eq2} and~\eqref{Aug29-eq3}-\eqref{Dec20-eq3}; the gradient field is indicated by vectors (red); and the limit cycle is in black.}\label{vdP1}
\end{figure}

The region between the
isolines $V_\ell(x)=\gamma_\ell$ and $V_u(x)=\gamma_u$ is the positively invariant set $\mathcal{B}$ shown in
Figure~\ref{vdP1}. Since $\mathcal{B}$ does not contain fixed points, we can apply, as an alternative to Li\'enard's theorem, the Poincar\'e-Bendixson Criterion. This implies that $\mathcal{B}$ must contain an asymptotically stable limit cycle. The limit cycle is shown in Figure~\ref{vdP1}. Note that the criterion only holds in the case of a second-order system. In the next section, we present a result for systems of any dimension that guarantees existence of exponentially stable limit cycles in a bounded region when the latter does not contain fixed points
and prove exponential stability of the van der Pol limit cycle.

\section{Exponential stability of limit cycles}\label{exponential}

In this section, we present dynamical systems with a certain contraction property. This property leads to stability of limit
sets. In~\cite{DCL49}, D.~C.~Lewis studied autonomous dynamical systems with a certain local contraction property. The idea behind it is that if trajectories remain in a bounded region and the distance between any two decreases with time then there exists a unique exponentially stable equilibrium point in that region. Bo T.~Stenstr\"om in~\cite{BS62} and Peter Giesl in~\cite{Giesl204} then provided sufficient conditions for the existence of a unique and
exponentially stable periodic orbit (see also~\cite{Borg60,Hartman}). Their concept is similar to the one by D.~C.~Lewis, with an important extension, namely, relaxing the requirement of contraction in the direction of the trajectory (see condition (\ref{positivity}) in Theorem~\ref{giesl}).

Consider the following system
\begin{equation}\label{may10-1}
\dot x=f(x),\ x\in\mathcal{B}\subseteq\mathbb{R}^n
\end{equation}
where $\mathcal{B}$ is a compact, connected and positively
invariant set of (\ref{may10-1}). First, we provide the definition of a \emph{Riemannian metric} for all $x\in\mathcal{B}$.

\begin{definition}\label{giesldef}{\bf (Definition 4
of~\cite{Giesl204})} A matrix-valued $C^1$-function $M(x)$ will be called a \emph{Riemannian metric} if $M(x)$ is a symmetric and positive definite matrix for all $x\in\mathcal{B}$.
\end{definition}

In this paper, whenever we require $M(x)>0$, we assume that it is a Riemannian metric. In the following, we present an important result by Giesl on the exponential stability of the periodic solution of autonomous dynamical systems, where $\mathcal{B}$ is as previously defined.

\begin{theorem}\label{giesl}{\bf(Theorem 5
of~\cite{Giesl204})} Let $M(x)$ be a Riemannian metric with its
directional derivative $M'(x)$ given by:
\begin{equation}
M'(x)_{(i,j)}=\sum_{k=1}^n\frac{\partial M(x)_{(i,j)}}{\partial
x_k}f_k(x). \nonumber
\end{equation}
If for all $x\in\mathcal{B}$,
\begin{equation}
w^T\left(M(x)\frac{\partial f(x)}{\partial
x}+\frac{1}{2}M'(x)\right)w<0,\ \forall w\ \text{\rm{such that} }
w^TM(x)f(x)=0\ \label{positivity}
\end{equation}
then (\ref{may10-1}) has a unique, exponentially asymptotically stable
periodic orbit in $\mathcal{B}$.
\end{theorem}

The theory of dynamical systems with a certain local contraction property provides an extension of Lyapunov stability theory; for example, in~\cite{AUGUST2011}, contraction theory was used to establish
global complete
synchronisation of coupled identical oscillators. Moreover, through \verb"SOSTOOLS" it offers means to computationally perform a search for conditions that guarantee the existence of an exponentially stable limit cycle, as we show in the following. Importantly, by using the theory developed in Section \ref{finding}, we restrict the search space to a positively invariant set of the state space (as opposed to its entirety), thus facilitating the search for the relevant SOS polynomials.

\subsection{Proving exponential stability of limit cycles using SOS decompositions}

If the vector field of (\ref{may10-1}) is polynomial, then (\ref{positivity}) is a conditional positivity condition on multivariate polynomials. It is well known that such conditions are difficult to check in general. However,
as shown above, by replacing the positivity conditions with SOS conditions,
we can use \verb"SOSTOOLS" to
provide sufficient conditions (`certificates') for the existence of a limit cycle in $\mathcal{B}$. In order to use the SOS framework, we first reformulate condition (\ref{positivity}).
\begin{theorem}
If there exists a (not necessarily) SOS polynomial $\alpha(x)$ such that
\begin{equation}\label{alpha1}
w^T\left(M(x)\frac{\partial f(x)}{\partial x}+\frac{1}{2}M'(x)\right)w - \alpha(x)\left(w^T M(x)f(x)\right)^2<0,\
\forall w\in\mathbb{R}^n,\ \forall x\in\mathcal{B}
\end{equation}
then~(\ref{positivity}) holds.
\end{theorem}
\begin{proof}
If (\ref{alpha1}) holds then for all $w$, such that $w^T M(x)f(x)=0$, (\ref{positivity}) holds.
Note that %the
polynomial function $\alpha(x)$ is only an aid to the computation, when $w^T M(x)f(x)\neq0$.
\end{proof}

As shown in the PhD dissertation of the first author \cite{August07}, if we define $v=M(x)w$ then~(\ref{alpha1}) becomes:
\begin{equation}\label{alpha2}
v^T\left(\frac{\partial f(x)}{\partial x}M(x)^{-1}-\frac{1}{2}\left(M(x)^{-1}\right)' - \alpha(x)f(x) f(x)^T\right)v<0,\
\forall v\in\mathbb{R}^n,\ \forall x\in\mathcal{B}
\end{equation}
where we have used the fact that
\begin{eqnarray*}
(M(x)M(x)^{-1})'=0=M'(x)M(x)^{-1}+M(x)(M(x)^{-1})'\\ \Longrightarrow
-(M(x)^{-1})'=M(x)^{-1}M'(x)M(x)^{-1}
\end{eqnarray*}
and we note that $M(x)>0$ if and only if $M(x)^{-1}>0$.

Let us consider a positively invariant set constructed using the approach presented in Section~\ref{finding}, where  $\mathcal{B}=\{x\in \mathbb{R}^n\ |\ B_I(x)\geq0\ \mathrm{and}\ B_O(x)\leq0\}$,
$B_I(x)=V_\ell(x)-\gamma_\ell$, $B_O(x)=V_u(x)-\gamma_u$, and both can be computed as above.
If $f(x)$, $B_I(x)$ and $B_O(x)$ are given by polynomial functions, we can use \verb"SOSTOOLS" to search for polynomial functions $M(x)$ and $\alpha(x)$ that fulfil SOS conditions derived from~(\ref{alpha2}). This is a feasibility problem that can be implemented as
\begin{eqnarray}
\mathrm{given}\quad && f(x),\ B_I(x),\ B_O(x),\ \delta>0\nonumber\\
\mathrm{search\ for}\quad && N(x),\ \alpha(x),\ p_1(v),\ p_2(v)\nonumber \\
\mathrm{subject\ to}\quad && p_1(v),\ p_2(v),\
v^\mathrm{T}\left(N(x)-\delta I\right)v\ \mathrm{are\ SOS\ for\ all}\ v,x\nonumber\\
&& -v^T\left(\frac{\partial f(x)}{\partial x}N(x)-\frac{1}{2}N'(x) -
\alpha(x)f(x) f(x)^T-\delta I\right)v +  \nonumber\\
&&\qquad\qquad \qquad -p_1(v)B_I(x)+p_2(v)B_O(x)\ \mathrm{is\ SOS\ for\ all}\ v,x\nonumber\\
\label{eq:feasibility_Giesl}
\end{eqnarray}
where $N(x)=M^{-1}(x)$. If a solution is found in terms of SOS decompositions then this provides a sufficient condition for the existence of an exponentially stable limit cycle in $\mathcal{B}$. Such a certificate is therefore a proof of the existence of the limit cycle. Clearly, these are not necessary conditions and the absence
of a solution says nothing about the existence of the limit cycle.
We now illustrate the method establishing the existence and exponential stability of the limit cycle for the van der Pol oscillator.

\subsection{Application to the van der pol oscillator: Proving exponential stability of its limit cycle}

Recall the set of equations for van der Pol oscillator given in~\eqref{van der Pol}.
Contraction analysis and \verb"SOSTOOLS" were used to establish global exponential stability of the origin for $k>0$ in~\cite{ParilloSlotine08}. (The origin is exponential stable if (\ref{positivity}) hold for all $w$ not only for those for which $w^TM(x)f(x)=0$.)
As discussed in Section~\ref{sec:vanderPol_LC},
for $k<0$, the existence and stability of the periodic orbit can be proved using Li\'enard's theorem.
As shown above, we can obtain such a region $\mathcal{B}$ as the region between
$B_I(x)=V_\ell-\gamma_\ell=0$ and $B_O(x)=V_u-\gamma_u=0$, where
\begin{eqnarray}
B_I(x)&=&-4-4.281x_1x_2+4.93x_1^2+2.213x_2^2-1.578x_1^4+0.012x_2^4+0.202x_1^6\nonumber\\
&&+1.07x_1^3x_2-0.048x_1^2x_2^2-0.483x_1x_2^3\label{Aug29-eq6}
\end{eqnarray}
and
\begin{eqnarray}
B_O(x)&=&-41+8.361x_1^2-11.679x_1x_2+4.925x_2^2-4.098x_1^4+3.240x_1^3x_2\nonumber\\
&&+0.710x_1^2x_2^2-0.063x_1x_2^3+0.050x_2^4+0.781x_1^6+0.298x_1^4x_2^2\nonumber\\
&&-0.434x_1^3x_2^3+0.283x_1^2x_2^4-0.089x_1x_2^5+0.012x_2^6.\label{Aug29-eq5}
\end{eqnarray}
Note that while we used $B_O(x)$ as defined by (\ref{Aug29-eq2}) and (\ref{Dec20-eq2}), for computational reasons, we chose a simpler description for $B_I(x)$ given by (\ref{Aug29-eq6}) instead of using our tighter expressions (\ref{Aug29-eq3}) and (\ref{Dec20-eq3}). We now use the boundaries~\eqref{Aug29-eq6} and~\eqref{Aug29-eq5} to establish the existence and exponential stability of a limit cycle of~(\ref{van der Pol}) for $k<0$.

\begin{figure}[h]
\begin{picture}(0,200)\hspace{-4cm}
\includegraphics[width=0.75\textwidth]{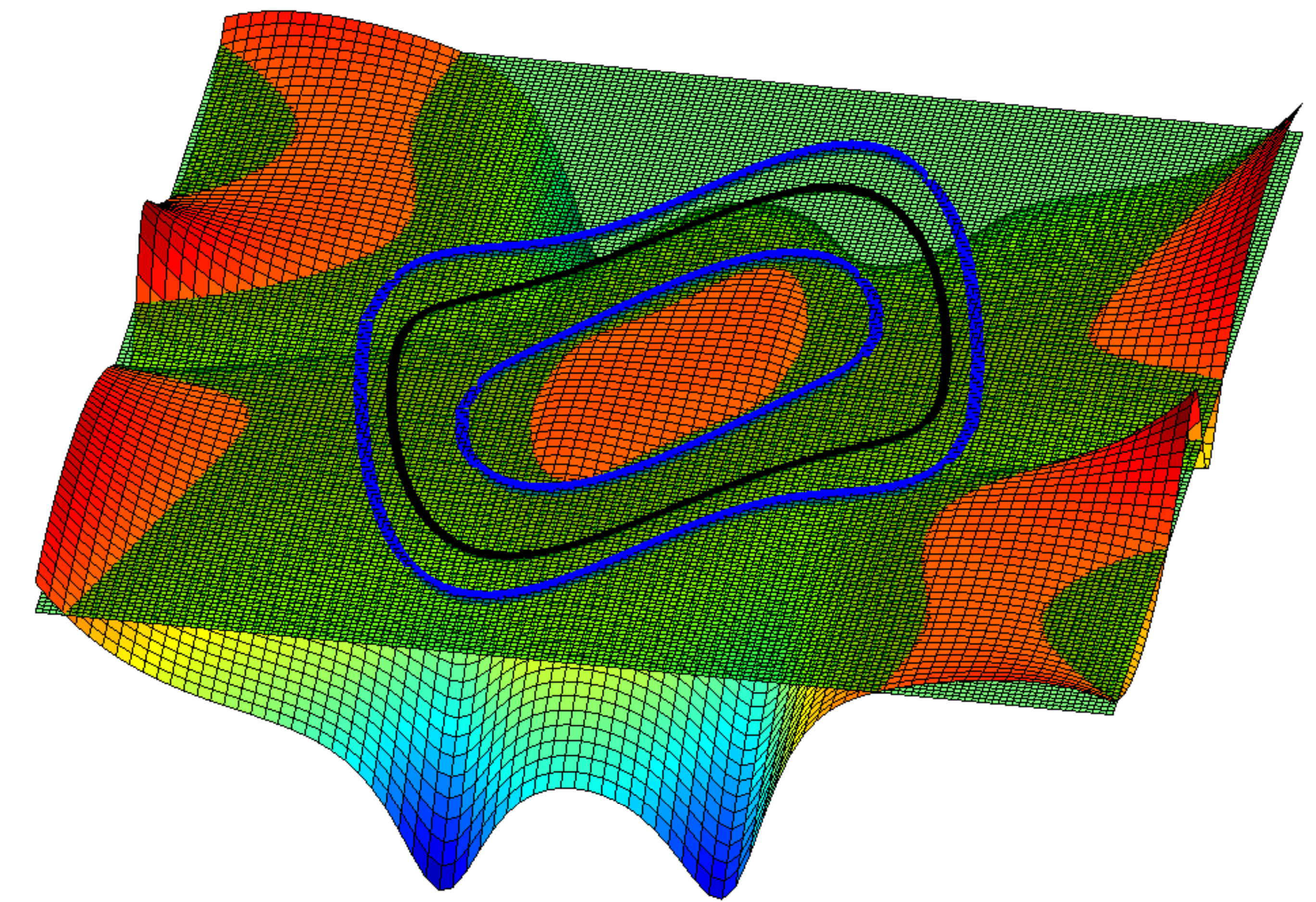}
\end{picture}
\caption{\textbf{Exponential stability of the limit cycle of the van der Pol oscillator.} The limit cycle is shown in black. The set $\mathcal{B}$ is the area between the blue curves. Here,
(\ref{alpha2}) holds if $L_M(x)<0$, where $L_M(x)=\frac{\partial f(x)}{\partial x}N(x)-\frac{1}{2}N'(x) - \alpha(x)f(x) f(x)^T$. The contour surface shown corresponds to the largest eigenvalue of $L_M(x)$, and the plane corresponds to the zero level.}\label{vdP2}
\end{figure}

Using \verb"SOSTOOLS", we solve the feasibility problem~\eqref{eq:feasibility_Giesl} to obtain polynomial functions $M(x)$, and
$\alpha(x)$. For $k=-1$ and $\delta=0.0001$, we
obtain:
\begin{eqnarray*}
N_{11}(x)&=& 2.759+1.080x_1^2+0.489x_1x_2+0.032x_2^2+0.514x_1^4+0.021x_1^3x_2\\
&&+0.002x_1^2x_2^2-0.006x_1x_2^3+0.001x_2^4\\
N_{12}(x)&=&-0.746-0.933x_1^2+1.652x_1x_2+0.033x_2^2-0.266x_1^4+0.014x_1^3x_2\\
&&+0.015x_1^2x_2^2-0.0002x_1x_2^3-0.001x_2^4\\
N_{22}(x)&=& 0.517+0.593x_1^2+0.717x_1x_2+3.900x_2^2+1.929x_1^4-0.184x_1^3x_2\\
&&+0.304x_1^2x_2^2+0.139x_1x_2^3+0.028x_2^4\\
\alpha(x)&=&
-1.592-2.310x_1^2-0.117x_1x_2-0.396x_2^2-0.977x_1^4-0.105x_1^3x_2\\
&&-0.488x_1^2x_2^2-0.151x_1x_2^3-0.296x_2^4.
\end{eqnarray*}
This result guarantees the existence of an exponentially stable
limit cycle for the van der Pol oscillator in $\mathcal{B}$ (Figure~\ref{vdP2}). The computation time is on the scale of a few
seconds on a standard PC. We have obtained similar results for several values of $k<0$.

%\newpage

\subsection{Application to three-dimensional dynamical systems}
In the following we apply our method to three-dimensional dynamical systems. This is important, since standard analytical methods to prove existence and stability of a limit cycle are constrained to two-dimensional ones. Here, we show that by using the theory developed by Peter Giesl in conjunction with our computational implementation we can search for a limit cycle and prove its exponential stability also in higher dimensional space.
\subsubsection{3D example 1}
\begin{figure}[h]
%\begin{picture}(0,260)\hspace{-1.6cm}
\begin{center}
\includegraphics[width=1\textwidth]{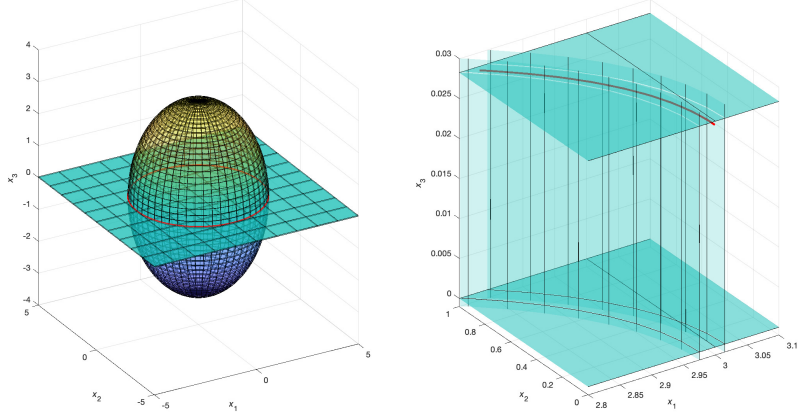}
\end{center}
%\end{picture}\vspace{-.5cm}
\caption{\textbf{Exponential stability of the limit cycle of the system (\ref{GieslDec14}).} See the text for our approach to prove existence and exponential stability of the limit cycle of (\ref{GieslDec14}). The numerical simulation of the limit cycle is shown in red.
We also show the bounding surfaces of the positively invariant set $\mathcal{B}_1$~\eqref{eq:Giesl_invariant_set}.
The right panel provides a close-up.}
\label{Orbits3D}
\end{figure}

Consider the following system from~\cite{Giesl03}:
\begin{eqnarray}\label{GieslDec14}
\dot x_1&=&(9-x_1^2-x_2^2-x_3^2)x_1-x_2\nonumber\\
\dot x_2&=&(9-x_1^2-x_2^2-x_3^2)x_2+x_1+0.25\nonumber\\
\dot x_3&=&(-x_1^2-x_2^2-x_3^2)x_3+0.25
\end{eqnarray}
System (\ref{GieslDec14}) has a unique real equilibrium point $(-0.0033,-0.0287,0.6295)$, which is unstable. It can be shown that if  $x_1^2(0)+x_2^2(0)+x_3^2(0)\geq\beta$ and $0\leq x_3(0)\leq \frac{0.25}{\beta}$, then $0\leq x_3(t)\leq \frac{0.25}{\beta}$ for all $t>0$. Using this fact, we can use \verb"SOSTOOLS" and the approach described previously to show: (i) that the set
\begin{align}
\label{eq:Giesl_invariant_set}
\mathcal{B}_1= \left \{x\in \mathbb{R}^3\ |\ 8.85 \leq x_1^2+x_2^2+x_3^2 \leq 9.09, \,
0\leq x_3\leq \frac{0.25}{8.85} \right \}
\end{align}
is positively invariant for (\ref{GieslDec14}); and, (ii) that there exists a matrix $N(x)>0$ and SOS functions $p_1(x)$, $p_2(x)$,
and $p_3(x)$ such that the following holds:
\begin{align}\label{Dec16eq1}
-J(x)N(x) + p_3f(x)f(x)^\mathrm{T}+p_1(x)(x_1^2+x_2^2+x_3^2-8.85)(x_1^2+x_2^2+x_3^2-9.09) \, I+
\nonumber\\
 +p_2(x) \, x_3\left(x_3-\frac{0.25}{8.85}\right)I-0.001\, I \geq0,
\end{align}
where $I$ is the identity matrix, $f(x)=[\dot x_1,\ \dot x_2,\ \dot x_3]^\mathrm{T}$, and $J(x)=\frac{\partial f(x)}{\partial x}$ denotes the Jacobian of $f(x)$.
We find that~(\ref{Dec16eq1}) is fulfilled
with
\begin{eqnarray*}
N(x) &=& \mathrm{diag}([0.662,\ 0.6612,\ 1.973])\\
p_1(x) &=& 0.108+0.0851x_3^2+0.0449x_2^2-0.0029x_1+0.0452x_1^2\\
p_2(x) &=&  0.5605+0.1135x_3+3.4503x_3^2+0.0122x_2+0.0002x_2x_3+3.8665x_2^2\\
&&-0.0481x_1-0.0005x_3x_1-0.002x_2x_1+3.8726x_1^2\\
p_3(x)&=&0.0994
\end{eqnarray*}
Therefore, this implies that there exists an exponentially stable limit cycle in $\mathcal{B}_1$, as shown in Figure~\ref{Orbits3D}.  Note that in this case it is enough to choose $N(x)$ as a diagonal matrix of constants, thus reducing the computational cost of the optimisation.

\subsubsection{3D example 2}
Consider the following system, which has a more complex limit cycle than the system~\eqref{GieslDec14}:
\begin{eqnarray}\label{3DDec14}
\dot x_1&=&x_1-x_2-x_1(x_1^2+x_2^2+x_3^2)\nonumber\\
\dot x_2&=&x_1+x_2-x_2(x_1^2+x_2^2-x_3^2)\nonumber\\
\dot x_3&=&x_2^2-x_3
\end{eqnarray}
This system has a unique real equilibrium point at the origin, which is unstable. Note that if $0\leq x_3(0)$ initially then $0\leq x_3(t)$ for all $t>0$.
Moreover, the cylinder given by $V(x)=1$ and $x_3\in\mathbb{R}$, where $V(x)=x_1^2+x_2^2$, is positively invariant, since
$$\dot V(x)=(x_1^2+x_2^2)(1-x_1^2+x_2^2)=V(x)(1-V(x)).$$

\begin{figure}[h]
\begin{center}
\includegraphics[width=.9\textwidth]{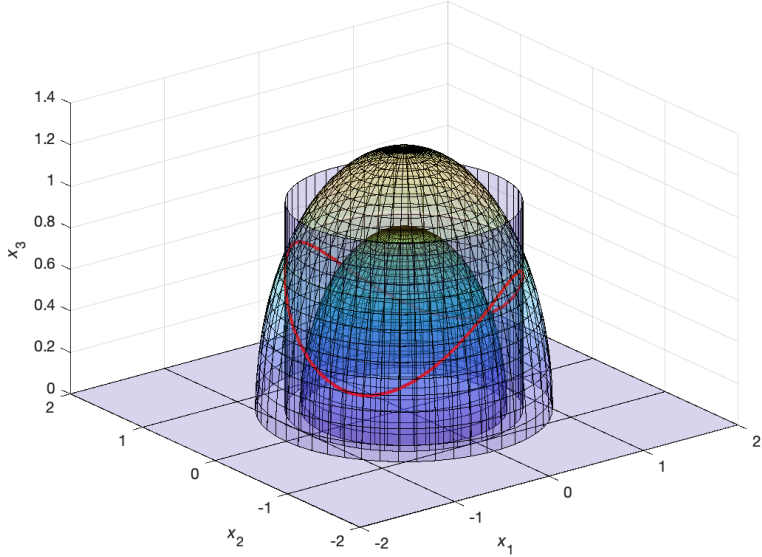}
\end{center}
\caption{\textbf{Exponentially stable limit cycle of the system~(\ref{3DDec14}).} See text for our approach to prove existence and exponential stability of the limit cycle. The numerical simulation of the limit cycle is shown in red.
We also show the bounding surfaces of the set $\mathcal{B}_2$ given by~\eqref{eq:pos_inv_set_3D_ex2}, i.e.,
the cylinder defined by $x_1^2+x_2^2=1$,
the half-spheres defined by $x_1^2+x_2^2+x_3^2=\delta$ with $\delta=\{0.75,1.56\}$, and the plane $x_3 = 0$.}
\label{Orbits3D2}
\end{figure}

Using this fact, we can use \verb"SOSTOOLS" and the approach described previously to show, first, that the set
\begin{align}
\label{eq:pos_inv_set_3D_ex2}
    \mathcal{B}_2 = \left \{x\in \mathbb{R}^3\ |\ 0.75 \leq x_1^2+x_2^2+x_3^2 \leq 1.56, \,\, x_1^2+x_2^2 \leq 1, \,\,  x_3 \geq 0
   \right \}
\end{align}
is positively invariant for (\ref{3DDec14}) and, then, that there exists a matrix $N(x)>0$ and SOS functions $p_1(x)$, $p_2(x)$,
and $p_3(x)$ such that the following holds,
\begin{eqnarray}\label{Dec16eq2}
&& \frac{1}{2}N^\prime(x)-J(x)N(x)
+p_1(x)(x_1^2+x_2^2+x_3^2-0.75)(x_1^2+x_2^2+x_3^2-1.56) \, I+ \nonumber\\
&& \qquad -p_2(x) x_3 \, I+ p_3(x) f(x)f(x)^\mathrm{T}+p_4(x)(x_1^2+x_2^2-1) \, I-0.001 \, I \geq0.
\end{eqnarray}
Then, (\ref{Dec16eq2}) implies that there exists an exponentially stable limit cycle in $\mathcal{B}_2$ (Figure~\ref{Orbits3D2}).
For completeness, we also report here the polynomials obtained for this system:
\begin{eqnarray*}
N(x) &=& \mathrm{diag}([0.4283,\ 0.5122,\ 2.7938])\\
p_1(x) &=& 0.403x_1^2 + 0.078x_1x_2  + 0.406x_2^2 + 0.325x_3^2 + 0.443\\
p_2(x)   &=&  0.8589\\
p_3(x) &=& 0.1168\\
p_4(x)&=& -0.910x_1^2 + 0.232x_1x_2 - 1.196x_2^2 + 1.527x_3^2 + 0.824
\end{eqnarray*}

\section{Conclusion}\label{conc}

In this paper, we have provided a computational means to search for positively invariant sets of polynomial dynamical systems using \verb"SOSTOOLS". We have exemplified the search through applications to different physical and biological systems of significance. For the Lorenz system, we obtained an attractive set and compared to a classic set from the literature, and combined both to obtain a smaller region that bounds the solution trajectories of the system. For the FitzHugh-Nagumo model of neuronal dynamics, we obtained positively invariant sets and used them to define a necessary condition for switching behaviour and to bound the maximal extension of the excursion after excitation.
Whilst a polynomial vector field that is globally asymptotically stable but does not have an analytic Lyapunov function was presented
in~\cite{Ahmadi18}, it is possible that there exist a Lyapunov function that is a SOS function and defines a positively invariant set for such vector fields. We leave this question as an open direction of future research.

Obtaining a positively invariant set for the van der Pol oscillator that did not contain an equilibrium point implied the existence of a stable limit cycle. Importantly, defining this set was also necessary to prove exponential stability of the limit cycle using further results presented in this paper.
In particular, we computationally implemented previously developed theoretical results by Giesl~\cite{Giesl03,Giesl204}  and used \verb"SOSTOOLS" to search for SOS decompositions that provide sufficient conditions for the existence of exponentially stable periodic solutions of polynomial dynamical systems. We applied the method to the van der Pol system, as well as to two higher dimensional systems (of dimension 3). This approach thus provides a computational means to search for limit cycles and prove their stability for systems of arbitrary dimension. Note, in contrast, that
references~\cite{ParilloSlotine08} and~\cite{MANCHESTER201432} use a particular 2-dimensional system when applying their approach, possibly because searching for stability certificates for higher order systems is hampered if the space of the search is not restricted. We remark that our approach
could be extended to deal with limit cycles that are not centred around the origin by first applying the method presented in~\cite{SCHLOSSER2021109900} to locate the limit cycle and then applying a coordinate shift. Finally, we would like to note that the approach presented here can also be extended to systems whose dynamics are described by rational functions by reformulating the problem (informally, by multiplying the equations by all denominators) at an increased computational cost.

\section*{Acknowledgments} Research funded by EPSRC. MB acknowledges funding from EPSRC grant EP/N014529/1 supporting the EPSRC Centre for Mathematics of Precision Healthcare. The authors would like to thank Pablo Parrilo and the anonymous reviewers for helpful comments about the paper.

% You may incorporate your references as follows in your main tex file.
% Using BibTex is not recommended but can be handled.

\medskip
% The data information below will be filled by AIMS editorial staff
Received October 2021; revised July 2022; early access August 2022.
\medskip

\end{document}